\documentclass[12pt]{article}

\usepackage[T1]{fontenc}
\usepackage[latin1]{inputenc}

\usepackage{booktabs}
\usepackage{amssymb,latexsym}
\usepackage{a4}
\usepackage[all]{xy}
\usepackage{color}

\usepackage[normalem]{ulem}
\usepackage{bbm}

\usepackage{rotfloat}

\begin{document}
\newfont{\blb}{msbm10 scaled\magstep1}

\newtheorem{theo}{Theorem}

\newtheorem{defi}[theo]{Definition}
\newtheorem{prop}[theo]{Proposition}
\newtheorem{lemm}[theo]{Lemma}
\newtheorem{coro}[theo]{Corollary}
\pagestyle{myheadings}
\date{}
\author{Ren{\'e} Hartung \\ 
Mathematisches Institut\\
Georg-August Universit\"{a}t zu G\"{o}ttingen \\
37073 G\"{o}ttingen, Germany \\
email: rhartung@uni-math.gwdg.de \\
\mbox{}\\
Gunnar Traustason \\
Department of  Mathematical Sciences, \\
University of Bath, \\
Bath BA2 7AY,
UK\\
email:gt223@bath.ac.uk}

\title{Refined solvable presentations for polycyclic groups}
\maketitle
\begin{abstract}
%
We describe a new type of polycyclic presentations, that we will
call refined solvable presentations, for polycyclic groups. These
presentations are obtained by refining a series of normal subgroups with
abelian sections. These presentations can be described effectively by 
presentation maps which yield the basis data structure to define a polycyclic 
group in computer-algebra-systems like {\scshape Gap} or 
{\scshape Magma}. We study refined solvable
presentations and, in particular, we obtain consistency criteria for them. This 
consistency  implementation demonstrates that it is often faster than the 
existing methods for polycyclic groups. 
%
%
\end{abstract}

%
\section{Introduction} 
A group $G$ is polycyclic if there exists a finite series of subnormal subgroups
$G = G_1 \unrhd G_2 \unrhd \ldots \unrhd G_m \unrhd G_{m+1} = \{1\}$ so that 
each section $G_i / G_{i+1}$ is cyclic. Polycyclic groups play an important role
in group theory as, for instance, each finite group with odd order is
polycyclic. Moreover, polycyclic groups form a special class of finitely
presented groups for which various algorithmic problems are solvable.
For instance, it is well-known that the word problem in a polycyclic
group is solvable. More precisely, a polycyclic group $G$ can be described by
a polycyclic presentation. This is a finite presentation with generators
$\{a_1,\ldots,a_m\}$ and relations of the form
$$
  \begin{array}{rcll}
  a_i^{r_i} &=& a_{i+1}^{\alpha_{i,i+1}} \cdots a_m^{\alpha_{i,m}},& i\in {\mathcal I}\\
  a_i^{-1} a_j a_i &=& a_{i+1}^{\beta_{i,j,i+1}} \cdots a_m^{\beta_{i,j,m}},&1\leq i<j\leq m\\
  a_i^{-1} a_j^{-1} a_i &=& a_{i+1}^{\gamma_{i,j,i+1}} \cdots a_m^{\gamma_{i,j,m}},&1\leq i<j\leq m,\: j\not\in {\mathcal I}\\
  a_i a_j a_i^{-1} &=& a_{i+1}^{\delta_{i,j,i+1}} \cdots a_m^{\delta_{i,j,m}},&1\leq i<j\leq m,\: i\not\in {\mathcal I}\\
  a_i a_j^{-1} a_i^{-1} &=& a_{i+1}^{\varepsilon_{i,j,i+1}} \cdots a_m^{\varepsilon_{i,j,m}},&1\leq i<j\leq m,\: i,j\not\in {\mathcal I}
  \end{array} $$
for a subset ${\mathcal I} \subseteq \{1,\ldots,m\}$ and integers $\alpha_{i,\ell},
\beta_{i,j,\ell}, \gamma_{i,j,\ell}, \delta_{i,j,\ell},\varepsilon_{i,j,\ell}\in {\mathbb
Z}$ that satisfy $0 \leq \alpha_{i,\ell}, \beta_{i,j,\ell},
\gamma_{i,j,\ell}, \delta_{i,j,\ell},\varepsilon_{i,j,\ell} < r_i$ whenever
$\ell\in {\mathcal I}$ holds.  For further details on polycyclic presentations we
refer to Section 9.4 of [14]. \\ \\
Given any finite presentation of a polycyclic group, the polycyclic
quotient algorithm [11,12] allows one to compute a
polycyclic presentation defining the same group. If, additionally,
the polycyclic group is nilpotent, than any finite presentation can be
transformed into a polycyclic presentation with the nilpotent quotient
algorithm [13]. We further note that even certain infinite
presentations (so-called finite $L$-presentations; see [2])
of a nilpotent and polycyclic group can be transformed into a polycyclic
presentation [3]. We may therefore always assume that a polycyclic
group is given by a polycyclic presentation. \\ \\
In the group $G$, every element is represented by a word $a_1^{e_1}
a_2^{e_2} \cdots a_m^{e_m}$ with $0\leq e_i < r_i$ whenever $i \in
{\mathcal I}$ holds. If this representation is unique, then the polycyclic
presentation is consistent and it yields a normal form for elements in
the group. This is a basis for symbolic computations within polycyclic
groups. Various strategies for computing normal forms in a polycyclic
group have been studied so far [10,16,6,1]. The current
state of the art algorithm is \emph{collection from the left}. But it is
known that even `collection from the left' is exponential in the number
of generators [10]; see also [1]. \\ \\
In this paper, we concentrate on refined solvable presentations as a special
class of polycyclic presentations that we describe in
Section 2. We
choose a finite series of normal subgroups so that
the sections are abelian. A refined solvable presentation will be a certain 
polycyclic presentation that refines this series.   Each weighted
nilpotent presentation, as used extensively in the nilpotent quotient
algorithms [13,3] and in [15], is of this type. A
solvable presentation can be described effectively by presentation
maps which we define in Section 2.  Presentation maps
can be considered as the basic data structure to define a polycyclic
group in computer-algebra-systems like {\scshape Gap} or {\scshape
Magma}. We obtain consistency criteria for refined solvable presentations in
Section 3. This consistency check has been implemented
in the {\scshape Nql}-package [8].  Our implementation shows
that the consistency checks for solvable presentations are often
faster than the general methods for polycyclic groups. As an example,
we consider nilpotent quotients of the Basilica group [7]
and the BSV group [4].  \\ \\
Fast algorithms for polycyclic groups are of special interest as, for
instance, the algorithm in [9] attempts to find periodicities
in the Dwyer quotients of the Schur multiplier of a group. In order
to observe these periodicities, the algorithm needs to compute with
polycyclic presentations with some hundreds of generators and therefore
fast algorithms for polycyclic groups are needed.

\section{Refined solvable presentations}
%
%
%
%
%
Let $G$ be a poly-cyclic group with a strictly ascending chain of
normal subgroups
       $$\{1\}=G_{0}<G_{1}<\cdots <G_{r}=G$$
where $G_{i}/G_{i-1}$ is abelian for $i=1,\ldots ,r$. Since each subgroup of a 
polycyclic group is finitely generated, we can choose a finite generating
set $X$ for $G$ which partitions as $X=X_{1}\cup X_{2}\cup \cdots \cup X_{r}$ 
such that 
        $$G_{i}/G_{i-1}=\bigoplus_{x\in X_{i}}\langle
        xG_{i-1}\rangle$$
for $i=1,\ldots ,r$ and where all the direct summands are non-trivial. We can 
furthermore make
our choice so that for each $x\in X_{i}$, either the order,
$o(xG_{i-1})$, of $xG_{i-1}$ is infinite or a power of a prime. Let
${\mathcal P}$ denote the set of all primes. For each $p\in {\mathcal P}$,
let
       $$X_{i}(p)=\{x\in X_{i}:\,o(xG_{i-1})\mbox{\ is a power of
       }p\}$$
and let
      $$X_{i}(\infty)=\{x\in X_{i}:\,o(xG_{i-1})=\infty\}.$$
Notice that the Sylow $p$-subgroup of $G_{i}/G_{i-1}$ is
      $$(G_{i}/G_{i-1})_{p}=\bigoplus_{x\in X_{i}(p)}\langle
      xG_{i-1}\rangle.$$
We order the generators in $X$ such that the generators in $X_{i}$
precede those in $X_{j}$ whenever $i<j$. Suppose that
$X=\{x_{1},\ldots ,x_{m}\}$ with $x_{1}<x_{2}<\ldots <x_{m}$.
 For each $x\in X_{i}$ let $n(x)=o(xG_{i-1})$. If $n(x)=\infty$,
let $\mbox{\blb Z}_{x}=\mbox{\blb Z}$ and otherwise let $\mbox{\blb
Z}_{x}=\{0,\ldots ,n(x)-1\}$. Each element $g\in G$ has a unique
normal form expression
     $$g=x_{m}^{r_{m}}x_{m-1}^{r_{m-1}}\cdots x_{1}^{r_{1}}$$
where $r_{i}\in \mbox{\blb Z}_{x_{i}}$. \\ \\
We next describe some relations that hold in the generators $x_{1},\ldots ,
x_{m}$. If $x\in X_{s}(p)$ then we get a {\it power relation} of the form
\begin{equation}
     x^{n(x)}=x_{m}^{\alpha_{x}(m)}\cdots x_{1}^{\alpha_{x}(1)}
\end{equation}
with $\alpha_{x}(i)\in \mbox{\blb Z}_{x_{i}}$ and where
$\alpha_{x}(i)=0$ if $x_{i}\not \in X_{1}\cup\cdots \cup X_{s-1}$. \\ \\
For each pair of generators $x,y\in X$ with $x<y$ we also get a
{\it conjugacy relation}
\begin{equation}
    x^{y}=x_{m}^{\beta_{(x,y)}(m)}\cdots
    x_{1}^{\beta_{(x,y)}(1)}
\end{equation}
where $\beta_{(x,y)}(i)\in \mbox{\blb Z}_{x_{i}}$. \\ \\
{\bf Remark}. There are three types of relations of the form (2). \\ \\
\underline{Type 1}. If $x,y\in
X_{s}$ then $x$ and $y$ commute modulo $G_{s-1}$ and thus we get
that $\beta_{(x,y)}(i)=0$ if $x_{i}\not\in X_{1}\cup \cdots \cup
X_{s-1}\cup \{x\}$ and that $\beta_{(x,y)}(i)=1$ if $x_{i}=x$. \\ \\
Now
suppose that $s<t$. \\ \\
\underline{Type 2}. If $x\in X_{s}(p)$ and $y\in X_{t}$ then
$x^{y}G_{s-1}\in (G_{s}/G_{s-1})_{p}$ and thus we get a relation of
the form (2) where $\beta_{(x,y)}(i)=0$ if $x_{i}\not \in X_{1}\cup
\cdots \cup X_{s-1}\cup X_{s}(p)$. \\ \\
\underline{Type 3}. Finally if $x\in X_{s}(\infty)$
and $y\in X_{t}$ then $x^{y}\in G_{s}$ and we get a relation of the
type (2) where $\beta_{(x,y)}(i)=0$ if $x_{i}\not\in
X_{1}\cup \cdots \cup X_{s}$. \\ \\
{\bf Remark}. By an easy induction on $m$, one can see that (1) and (2) also 
give us, for every pair of generators $x,y\in X$ such that $x<y$,  
a relation $x^{y^{-1}}=\mu(x,y)$, where $\mu(x,y)$ is a normal form 
expression. Thus using only relations (1) and the three types of relations
(2), we have a full information about $G$ and we can calculate inverses and 
products of elements of normal form and turn the result into a normal form
expression using for example collection from the left. \\ \\
The claim holds trivially for $m=1$. Now suppose that $m\geq 2$
and that the claim holds for all smaller values of $m$. Consider the subgroup
$H=\langle x_{1},\ldots ,x_{m-1}\rangle$. By the inductive hypothesis, every
element in $H$ can be turned into a normal form expression using only
relations (1) and (2). Now (2) gives us normal form expressions for 
$x_{1}^{x_{m}},\ldots ,x_{m-1}^{x_{m}}$ and this determines an automorphism
$\phi\in \mbox{Aut\,}(H)$ induced by the conjugation of $x_{m}$. This
then gives us $\phi^{-1}$ that gives us in turn normal form expressions
for $x_{1}^{x_{m}^{-1}},\ldots ,x_{m-1}^{x_{m}^{-1}}$. This finishes the proof
of the inductive step. \\ \\
The point about this is that the relations $x^{y^{-1}}=\mu(x,y)$ are not defining
relations but consequences of (1) and (2). So for a polycyclic group $G$ we
only need (1) and (2) to define it. For practical reasons we need however 
to determine the relations $x^{y^{-1}}=\mu(x,y)$ first to be able to perform 
calculations in $G$. At the end of section $3$, we describe an efficient method
for doing this for the polycyclic presentations that we are about to introduce
next, refined solvable presentations. \\ \\ \\
%
%
%
%
Suppose now conversely that we have a finite alphabet 
$X=\{x_{1},x_{2},\ldots ,x_{m}\}$ with an ordering $x_{1}< x_{2}<
\ldots < x_{m}$. Let $F$ be the free group on $X$. Partition $X$ into some 
disjoint  non-empty subsets
$X_{1},\ldots, X_{r}$ such that the elements of $X_{i}$ precede
those in $X_{j}$ whenever $i<j$. Then partition further each $X_{i}$
as a union of disjoint
subsets (most empty of course) \\ \\
      $$X_{i}=(\bigcup_{p\in  {\mathcal P}}X_{i}(p))\cup
      X_{i}(\infty).$$
Let $Y=Z\setminus \{x\in X:\,n(x)=\infty\}$ and $Z=\{(x,y)\in X\times X:\,
x<y\}$. We introduce three maps that we will refer to as {\it presentation
maps}. The first one is
             $$n:X\rightarrow \mbox{\blb N}\cup \{\infty\}$$
such that $n(x)=\infty$ if $x\in X_{i}(\infty)$ and $n(x)$ is a non-trivial
power of $p$ is $x\in X_{i}(p)$. The second presentation map is 
     $$\pi:Y\rightarrow F$$
where, if $x\in X_{s}(p)$, $\pi(x)=x_{m}^{\alpha_{x}(m)}\cdots x_{1}^{\alpha_{x}(1)}$ with
$\alpha_{x}(i)\in \mbox{\blb Z}_{x_{i}}$ and $\alpha_{x}(i)=0$
whenever $x_{i}\not\in X_{1}\cup\cdots \cup X_{s-1}$. Notice that these
are the conditions for the right hand side of the power relation (1). The
final presentation map is
     $$\delta:Z\rightarrow F$$
where $\delta(x,y)=x_{m}^{\beta_{(x,y)}(m)}\cdots x_{1}^{\beta_{(x,y)}(1)}$ and
the conditions for the right hand side of (2) above hold as indicated in the
remark that follows it. So we have a data that consists of an alphabet $X$
with a partition and three presentation maps. To this data we associate a
presentation with generators $x_{1},\ldots ,x_{m}$, power relations 
                 $$x^{n(x)}=\pi(x)$$
for any $x\in X$ such that $n(x)\not =\infty$, and conjugacy relations 
              $$x^{y}=\delta(x,y)$$
for each pair $(x,y)\in X\times X$ such that $x<y$. We call such a presentation
a {\it refined solvable presentation}. We have seen above that every polycyclic
group has a refined solvable presentation that is consistent. Conversely,
we are interested in criteria for a given refined solvable presentation 
to be a consistent presentation for a polycyclic group $G$. In other 
words we want the group $G$ to be polycyclic and we want every element $g\in G$
to have a unique normal form expression
           $$g=x_{m}^{r_{m}}\cdots x_{1}^{r_{1}}$$
with $r_{i}\in \mbox{\blb Z}_{x_{i}}$. In next section we describe such 
consistency criteria. \\ \\
{\bf Remark}. Notice that there are groups with a refined solvable 
presentation that are not polycyclic. Take for example two variables
$x_{1}<x_{2}$ and let $X_{1}=X_{1}(\infty)=\{x_{1}\}$, 
$X_{2}=X_{2}(\infty)=\{x_{2}\}$. Here $Y=\emptyset$ and $Z=\{(x_{1},x_{2})\}$.
For the presentation maps $n:X\rightarrow \mbox{\blb N}\cup \{\infty\}$ and
$\pi:Y\rightarrow F$, we must have $n(x_{1})=n(x_{2})=\infty$ and 
$\pi$ must be empty. Suppose we choose $\delta:Z\rightarrow F$ such that
$\delta(x_{1},x_{2})=x_{1}^{2}$. Then we get a presentation with two 
generators $x_{1},x_{2}$ and one relation
            $$x_{1}^{x_{2}}=x_{1}^{2}.$$
The resulting group is not polycyclic. The criteria that we will describe
in section 3 are thus not only consistency criteria but also criteria
for the resulting group to be polycyclic. 
\section{The consistency criteria}
Before establishing our consistency criteria, we first describe constructions
that are central to what follows. 
%
%
%
%
%
%
Suppose we have a polycyclic group $G=\langle X\rangle$ that has a consistent
refined solvable presentation as described above with a generating set
$X=\{x_{1},\ldots ,x_{m}\}$ that is partitioned as described in section 2 and
with presentation maps $n,\pi$ and $\delta$. Let $\phi\in \mbox{Aut\,}(G)$.
We will consider two situations where we can use this data to get a consistent
refined solvable presentation for a larger polycyclic group $\tilde{G}$. 
Add a new variable $x_{m+1}$ and extend our order on $\tilde{X}=X\cup
\{x_{m+1}\}$ such that $x_{m+1}$ is larger than the elements in $X$. Let
$\tilde{F}$ be the free group on $\tilde{X}$. Let $H$ be the semidirect
product of $G$ with a infinite cyclic group 
$C_{\infty}=\langle x\rangle$ where the action from $C_{\infty}$ on $G$
is given by $g^{x}=g^{\phi}$. \\ \\
For the first situation let $\tilde{G}=H$. We extend the presentation maps 
$n,\pi,\delta$ to $\tilde{n},\tilde{\pi},\tilde{\delta}$ so they involve 
$\tilde{X}$. We do this by letting $\tilde{n}(x_{m+1})=\infty$ and
        $$\tilde{\delta}(x_{i},x_{m+1})=x_{i}^{\phi}\mbox{\ \ (in a normal form
       expression in $x_{1},\ldots ,x_{m}$)}$$
for $i=1,\ldots ,x_{m}$. Notice that, since $n(x_{m+1})=\infty$, $\tilde{\pi}=
\pi$. The refined solvable presentation that we get using the extended
presentation maps has all the relations for $G$ together with $m$
extra relations
       $$x_{i}^{x_{m+1}}=\delta(x_{i},x_{m+1})=x_{i}^{\phi}$$
for $i=1,\ldots ,m$. A moments glance should convince the reader that this
is a refined solvable presentation for the polycyclic group $\tilde{G}=H$. \\ \\
{\bf Remark}. We haven't said anything above about the partition of 
$\tilde{X}=\{x_{1},\ldots ,x_{m+1}\}$.  The partition would be
into $\tilde{X}_{1}=X_{1},\ldots
,\tilde{X}_{r}=X_{r},\tilde{X}_{r+1}=\{x_{m+1}\}$. If furthermore
$x^{-1}x^{\phi}\in G_{r-1}$ for all $x\in X_{r}$ we could instead
choose a partition with $\tilde{X}_{1}=X_{1},\ldots
,\tilde{X}_{r-1}=X_{r-1},\tilde{X_{r}}=X_{r}\cup \{x_{m+1}\}$. \\ \\
The second situation is a variant of the first. 
Now suppose furthermore that for some integer $e\geq 2$, that
is a power of a prime $p$, and $g\in G$
we have that
\begin{eqnarray}
          a^{g} & = & a^{\phi^{e}}\mbox{\ \ (for all $a\in G$)} \\
           g^{\phi} & = & g
\end{eqnarray}
In this case $N=\langle g^{-1}x^{e}\rangle$ is a subgroup of the
centre of $H$. Let $\tilde{G}=H/N$. $G$ embeds naturally into
$\tilde{G}$ and we identify it with it's image. We now extend the presentation
maps $n,\pi,\delta$ to $\tilde{n},\tilde{\pi},\tilde{\delta}$ as follows.
First we let $\tilde{n}(x_{m+1})=e$ and $\tilde{\pi}(x_{m+1})$ be the
normal form expression for $g$ in $x_{1},\ldots ,x_{m}$. Finally as before
let $\tilde{\delta}(x_{i},x_{m+1})$ be the normal form expression of
$x_{i}^{\phi}$ in $x_{1},\ldots ,x_{m}$. The refined solvable presentation
with respect to the presentation maps $\tilde{n},\tilde{\pi}$ and 
$\tilde{\delta}$ is then a presentation with all the relations for 
$G$ and the extra relations
       $$x_{m+1}^{n(x_{m+1})}=\tilde{\pi}(x_{m+1})=g$$
together with 
       $$x_{i}^{x_{m+1}}=\tilde{\delta}(x_{i},x_{m+1})=x_{i}^{\phi}\mbox{\ \ (in a normal form
       expression in $x_{1},\ldots ,x_{m}$)}$$
for $1\leq i\leq m$.
Again it is clear that this is a refined solvable presentation for 
the polycyclic group $\tilde{G}=
H/N$. The remark above applies again for the partition in this case. \\ \\
\\
We now turn back to our task of finding a consistency criteria for
power-conjugate presentations of poly-cyclic groups. Suppose
$G=\langle x_{1},\ldots ,x_{m}\rangle$ is a poly-cyclic group with a
refined solvable presentation as described above. So we have some partition
of $X=\{x_{1},\ldots ,x_{m}\}$ and presentation maps $n,\pi,\delta$ giving
us relations
       $$x^{n(x)}=\underbrace{x_{m}^{\alpha_{x}(m)}\cdots
        x_{1}^{\alpha_{x}(1)}}_{\pi(x)}$$
for $x_{1}\leq x\leq x_{m}$ with $n(x)<\infty$ and 
       $$x^{y} = \underbrace{x_{m}^{\beta_{(x,y)}(m)}\cdots
                             x_{1}^{\beta_{(x,y)}(1)}}_{\delta(x,y)}$$
for $x_{1}\leq x<y\leq x_{m}$. For $k=0,1,\ldots ,m$, let $H_{k}$ be
the group satisfying the sub-presentation with generators
$x_{1},\ldots ,x_{k}$ and those of the relations involving only
$x_{1}\leq x<y\leq x_{k}$. The idea is to establish inductively
criteria for the refined solvable presentation for $H_{k}$ to be a consistent
presentation of a polycyclic group. The induction
basis $k=0$ doesn't need any work. Now suppose that we have already
obtained criteria for the refined solvable presentation for 
$H_{k}$, where
$0\leq k\leq m-1$, to be a consistent presentation of
a polycyclic group. Using the presentation map $\delta$ we define a function
$\delta(x_{k+1}):H_{k}\rightarrow H_{k}$ by first defining the
values of the generators as $x_{i}^{\delta(x_{k+1})}=\delta(x_{i},x_{k+1})$
for $i=1,\ldots ,k$. We then extend this to the whole of $H_{k}$ by
letting $\delta(x_{k+1})$ act on normal form expressions as follows
      $$(x_{k}^{r_{k}}\cdots x_{1}^{r_{1}})^{\delta(x_{k+1})}
     =(x_{k}^{\delta(x_{k+1})})^{r_{k}}\cdots
              (x_{1}^{\delta(x_{k+1})})^{r_{1}}.$$
Suppose the resulting map $\delta(x_{k+1})$ is an automorphism. If
$n(x_{k+1})=\infty$, we have that the presentation for $H_{k+1}$ is a
consistent presentation for the semidirect product of $H_{k}$ with the
infinite cyclic group $C_{\infty}=\langle x\rangle$ where 
$g^{x}=g^{\delta(x_{k+1})}$. Now suppose that $n(x_{k+1})\not
=\infty$. Using the second construction above and taking into
account conditions (3) and (4), we get a presentation for
$H_{k+1}$ that is a consistent presentation of a polycyclic group, provided that
\begin{eqnarray*}
      \pi(x_{k+1})^{\delta(x_{k+1})} & = & \pi(x_{k+1}) \\
      x_{i}^{\delta(x_{k+1})^{n(x_{k+1})}} & = & x_{i}^{\pi(x_{k+1})}
\end{eqnarray*}
for $i=1,\ldots ,k$. It remains to find criteria for
$\delta(x_{k+1})$ to be an automorphism. This problem we turn to
next. \\ \\ \\
Let $G=\langle X\rangle$ be a poly-cyclic group with a
consistent refined solvable presentation as described above. For $s=1,\ldots
,r$ let $G_{s}=\langle X_{1}\cup\cdots\cup X_{s}\rangle$,
$G_{s}(p)=\langle X_{1}\cup\cdots\cup X_{s-1}\cup X_{s}(p)\rangle$
and let $\tau(G_{s})=\langle X_{1}\cup\cdots
X_{s-1}\cup(\bigcup_{p\in{\mathcal P}}X_{s}(p))\rangle$. For each
$x\in X$ choose an element $x^{\phi}$ subject to the following
conditions:
\begin{eqnarray}
    x^{\phi}\in G_{i} & \mbox{if} & x\in X_{i} \\
    x^{\phi}\in G_{i}(p) & \mbox{if} & x\in X_{i}(p). \nonumber
\end{eqnarray}
We extend this to a map $\phi:G\rightarrow G$ by letting $\phi$ act
on normal form expressions as:
   $$(x_{m}^{r_{m}}\cdots x_{1}^{r_{1}})^{\phi}=
      (x_{m}^{\phi})^{r_{m}}\cdots (x_{1}^{\phi})^{r_{1}}.$$
Notice that the condition (5) implies that $\phi$ induces maps
$\phi_{s}:G_{s}\rightarrow G_{s}$, $s=1,\ldots ,r$,  where
$\phi_{s}=\phi|_{G_{s}}$. It also induces maps
$\phi_{(s,p)}:G_{s}(p)/G_{s-1}\rightarrow G_{s}(p)/G_{s-1}$ and maps
$\phi_{(s,\infty)}:G_{s}/\tau(G_{s})\rightarrow G_{s}/\tau(G_{s})$.
\\
\begin{lemm}
The map $\phi:G\rightarrow G$ is a homomorphism if and only if
\setcounter{equation}{0}
\begin{equation}
    \pi(x)^{\phi}= (x^{\phi})^{n(x)}\mbox{\ \ \ $(x_{1}\leq x\leq
    x_{m})$}
\end{equation}
and
\begin{equation}
 \mbox{}  x^{y \phi}=x^{\phi y^{\phi}}
    \mbox{\ \ \ \ \ $(x_{1}\leq x<y\leq x_{m})$}.
\end{equation}
$\phi$ is furthermore an automorphism if for $s=1,\ldots ,r$ we have
\begin{eqnarray}
      \mbox{det\,}(\phi_{(s,p)}) & \not = & 0\mbox{\ $($mod $p)$} \\
      \mbox{det\,}(\phi_{(s,\infty)}) & = & \pm 1. \nonumber
\end{eqnarray}
\end{lemm}
{\bf Proof}.\ \ Consider the homomorphism $\psi:F\rightarrow F$ on
the free group $F=\langle x_{1},\ldots ,x_{m}\rangle$ induced by the
values $x^{\psi}=x^{\phi}$ for $x_{1}\leq x\leq x_{m}$. Let $R$ be
the normal subgroup generated by the defining polycyclic relators for $G$.
This means that $G=F/R$. Then conditions (1) and (2) imply that
$R^{\psi}\leq R$ and thus $\psi$ induces a homomorphism on $G=F/R$.
This homomorphism is clearly the map $\phi$. \\ \\
The homomorphism $\phi$ is bijective if and only if the induced
linear maps $\phi_{(s,p)}$ and $\phi_{(s,\infty)}$ are bijective and
this happens if and only if condition (3) holds. $\Box$ \\ \\ \\
{\bf Remark}. The condition (1) in the lemma above is of course only relevant 
when $n(x)<\infty$. To avoid making the statement more complicated we can
decide that $\pi(x)=1$ and $u^{n(x)}=1$ for all $u\in G$ in the case
when $n(x)=\infty$. \\ \\ 
We now turn back again to the problem of establishing  
criteria for refined solvable presentations to be a consistent presentation 
of a polycyclic group. Let $G=\langle
x_{1},\ldots ,x_{m}\rangle$ be a group satisfying a
refined solvable presentation as described above with relations
\begin{eqnarray*}
        x^{n(x)} & = & \underbrace{x_{m}^{\alpha_{x}(m)}\cdots
        x_{1}^{\alpha_{x}(1)}}_{\pi(x)} \mbox{\ \ \ $(x_{1}\leq x\leq x_{m})$}\\
        x^{y} & = & \underbrace{x_{m}^{\beta_{(x,y)}(m)}\cdots
                             x_{1}^{\beta_{(x,y)}(1)}}_{\delta(x,y)}
                  \mbox{\ \ \ $(x_{1}\leq x<y\leq x_{m})$}.
\end{eqnarray*}
We let $H_{k}$ be the group satisfying the sub-presentation with
generators $x_{1},\ldots ,x_{k}$ and those of the relations where
$x_{1}\leq x<y\leq x_{k}$. We establish inductively 
criteria for the presentation for $H_{k}$ to be a consistent presentation of
a polycyclic group. Suppose this has been
achieved for some $k$. We want to add criteria so that the
presentation for $H_{k+1}$ is a consistent presentation for a polycyclic group.
We let
$\delta(x_{k+1}):H_{k}\rightarrow H_{k}$ be the map induced by the
values $x^{\delta(x_{k+1})}$ as described above. As we
pointed out, the presentation for $H_{k+1}$ is a consistent presentation 
of a polycyclic group if and only
if the map $\delta(x_{k+1})$ is an automorphism and that we have the
extra criteria that
\begin{eqnarray*}
      \pi(x_{k+1})^{\delta(x_{k+1})} & = & \pi(x_{k+1}) \\
      x_{i}^{\delta(x_{k+1})^{n(x_{k+1})}} & = &
      x_{i}^{\pi(x_{k+1})}.
\end{eqnarray*}
From Lemma 1 we have criteria for $\delta(x_{k+1})$ to be an
automorphism. Suppose that $x_{k+1}\in X_{s}$. Then
$\delta(x_{k+1})$ acts trivially on $G_{s}/G_{s-1}$ and so to
establish that $\delta(x_{k+1})$ is bijective we only need to show
that $\delta(x_{k+1})_{(t,p)}$ and
$\delta(x_{k+1})_{(t,\infty)}$ are bijective for $1\leq t<s$. \\ \\
For $z\in X$ let $r(z)$ be the integer such that $z\in X_{r(z)}$.
Adding up for $k=0,\ldots ,m-1$, we obtain the following consistency
criteria.
\begin{theo}
The refined solvable presentation for $G$ is a consistent presentation for
a polycyclic group if and only
if the following criteria hold. Firstly we must have for
all $x_{2}\leq z\leq x_{m}$ that
\setcounter{equation}{0}
\begin{eqnarray}
       \pi(z)^{\delta(z)} & = & \pi(z) \\
       \pi(x)^{\delta(z)} & = & (x^{\delta(z)})^{n(x)} \mbox{\ \ \
       \ \ \ \
       $(x_{1}\leq x<z)$} \\
       x^{\delta(z)^{n(z)}} & = & x^{\pi(z)} \mbox{\ \ \ \ \ \ \ \ \ \ \ \ \ $(x_{1}\leq
       x<z)$} \\
       \mbox{}    x^{y \delta(z)} & = & x^{\delta(z)y^{\delta(z)}}
        \mbox{\ \ \ \ \ $(x_{1}\leq x<y<z)$}.
\end{eqnarray}
We also need for $1\leq s<r(z)$ that
\begin{eqnarray}
      \mbox{det\,}(\delta(z)_{(s,p)}) & \not = & 0\mbox{\ $($mod $p$$)$} \\
      \mbox{det\,}(\delta(z)_{(s,\infty)}) & = & \pm 1. \nonumber
\end{eqnarray}

\end{theo}
\mbox{}\\
{\bf Remarks}. (1)  Recall that we established the consistency of the
polycyclic group 
$H_{k}$ recursively for $k=0,1,\ldots ,m$. So according to the proof
we should check (1)-(5) for $z=x_{2},\ldots ,x_{m}$ in ascending
order. If $z=x_{k+1}$ then the consistency of $H_{k+1}$ follows from
the consistency of $H_{k}$ together with relations (1)-(5) of Theorem
2 where $z=x_{k+1}$. So when doing the check for $z=x_{k+1}$ we
can assume that the presentation for $H_{k}$ is consistent. Using the
definition of $\delta(z)$ we first transform all the expressions in
(1)-(4) into expressions in $H_{k}$. Then we turn each side of the
equations into normal form in $H_{k}$ and compare. It is interesting
to note that (provided the check has been positive so far) $H_{k}$
has a consistent presentation and so the normal form in each case is
independent of how we calculate. We can however do the check in any order we like (and still
sticking to the assumption that $H_{k}$ has a consistent presentation).
The reason for this is that we will at some point reach the smallest
$z$ where the check fails (provided that we haven't got a negative
result in the mean time). Hence if the presentation is not a consistent
presentation of a polycyclic group, 
this will be recognised. \\ \\
(2) How does this approach compare to the existing ones.
%
Our approach is to consider functions $\delta(z)$ defined
on a group $G_{z}$ with a subpresentation (involving only the generators less 
than $z$). Modulo consistency of $G_{z}$ the conditions (1)-(5) in Theorem 2
are conditions for the map $\delta(z)$ to be an automorphism ((2), (4) and (5)) and 
for the resulting cyclic extension to have a consistent presentation ((1) and
(3)). The emphasis is thus on the function $\delta(z)$ rather than the group
operation (as in [14]). It is our belief that this viewpoint makes things
look a bit clearer. \\ \\
(3) It should be noted however that our conditions (1)-(4) have equivalent
criteria in the standard approach. See the list (*) in [14], page 424. The 'overlaps' (1),(2),(3) and (5) in that list correspond to (4),(2),(3) and (1) in 
Theorem 2. The condition (5) is however new and is a biproduct of working
with an ascending normal solvable series. In the standard approach one works
with a ascending subnormal series with cyclic factors. 
It should also be noted that the idea of obtaining consistency recursively for 
$H_{k}$, $k=0,\ldots, m$, through working with $\delta(z)$, is also implicit in [14] but is kept in the background within the proof. Our conditions (1)-(5) 
bring this to the surface.
%
\\ \\ \\
{\bf A method for obtaining inverse conjugation relations}. 
For practical checks using these consistency criteria one needs to
determine first normal form expressions $x^{z^{-1}}$ for $x<z<x_{m}$ 
(in order to be able to transform any expression in  $H_{k}$ to an
normal form expression). Note however that this is ofcourse only 
needed when $z$ is of infinite order. Another advantage of our approach is that it becomes
quite simple and effective to determine these after having produced all
the linear maps $\delta(z)_{(s,p)}$ and 
$\delta(z)_{(s,\infty)}$, $2\leq s\leq r$.
Suppose that $z\in X_{s}$ for some $2\leq s\leq r$. We now describe how to obtain
normal form expressions for
$x^{z^{-1}}$ recursively for $x<z$. \\ \\
We can suppose that we already know that the sub-presentation for the group $G^{*}$
generated by the generators $\{x\in X:\,x<z\}$ (using only the relations involving
these generators) is consistent. 
The presentation for $G^{*}$ is built around an ascending normal $z$-invariant series with
each factor either a finite abelian $p$-group or a finitely
generated torsion-free abelian group. \\ \\
Now suppose that we are looking at one such factor $K/H$ and that
the extra generators needed to generate $K$ are $y_{1},\ldots ,y_{e}$. We can suppose inductively
that
we have obtained normal form expressions for all $x^{z^{-1}}$ when
$x$ is a generator of $H$. We want to extend this
to $y_{i}^{z^{-1}}$ for $i=1,\ldots ,e$. \\ \\
Let $v_{1}=y_{1}H,\ldots ,v_{e}=y_{e}H$ be the generators of $K/H$.
Let $\phi$ be the automorphism on $K/H$ induced by the conjugation
action by $z$ and let $\psi$ be the inverse of $\phi$. Suppose
$\psi$ is represented by the matrix $B=(b_{ij})$. Since
$\phi(\psi(v_{i}))=v_{i}$, we have
      $$b_{ei}\phi(v_{e)}+\cdots +b_{2i}\phi(v_{2})+ b_{1i}\phi(v_{1})=v_{i}.$$
It follows that (using the presentation and calculating in $K$) we
get
    $$(y_{e}^{z})^{b_{ei}}\cdots (y_{2}^{z})^{b_{2i}}(y_{1}^{z})^{b_{1i}}=y_{i}u.$$
Where $u$ is a normal form expression in the generators of $H$ (and
we already know how $z^{-1}$ acts on $u$. It follows that
     $$y_{i}^{z^{-1}}=y_{e}^{b_{ei}}\cdots
     y_{2}^{b_{2i}}y_{1}^{b_{1i}}u^{-z^{-1}}.$$
\section{Implementation and some applications of our consistency checks}
We have implemented our consistency check in the {\scshape Nql} package [8]
of the computer-algebra-system {\scshape Gap}; see [5]. In
this section, we demonstrate how this method yields a significant speed-up
in checking consistency of large polycyclic presentations (with some hundreds of generators). For this
purpose, we consider nilpotent quotients of the Basilica group $\Delta$ from [7] and
the Brunner-Sidki-Vieira-Group $\mbox{BSV}$ from [4]. Both groups
are two-generated but infinitely presented. The Basilica group admits
the following infinite presentation
\[
  \Delta \cong \langle \{a,b\} \mid  [a,a^b]^{\sigma^{i}}, i\in\mbox{\blb N}_0\rangle
\]
where $\delta$ is the endomorphism of the free group over $a$ and $b$
induced by the mapping $a\mapsto b^2$ and $b\mapsto a$; see [7].
The $\mbox{BSV}$ group admits the infinite presentation
\[
  \mbox{BSV} \cong \langle \{a,b\} \mid [b,b^a]^{\varepsilon^{i}}, [b,b^{a^3}]^{\varepsilon^{i}}, i\in\mbox{N}_0 \rangle,
\]
where $\varepsilon$ is the endomorphism of the free group over $a$ and
$b$ induced by the mapping $a\mapsto a^2$ and $b\mapsto a^2b^{-1}a^2$.
The nilpotent quotient algorithm in [3] computes a 
weighted nilpotent presentation for the lower central series quotient
$G/\gamma_c(G)$ for a group $G$ given by an infinite presentation as
above (a so-called finite $L$-presentation; see [2]). A weighted nilpotent presentation is a polycyclic presentation
which refines the lower central series of the group. We note that the
weighted nilpotent presentations for the quotients $\Delta/\gamma_c\Delta$
and $\mbox{BSV}/\gamma_c\mbox{BSV}$ are refined solvable presentations. \\ \\
%
%
In order to verify consistency of a given polycyclic presentation,
the algorithm in [14, p. 424]  rewrites the overlaps of the rewriting rules
and compares the result; that is, the algorithm checks the underlying rewriting system for local confluence. As even the state of art algorithm 'collection from
the left' is exponential [10], the number of overlaps is a central bottleneck
here. There are improvements known which
make use of the structure of a polycyclic presentation in order to
reduce the number of overlaps. For instance, for weighted nilpotent
presentations, a weight function allows one to reduce the number of
overlaps significantly; see [14, p. 431]. \\ \\
Our method replaces some overlaps by the computation of determinants of
integral matrices and it can easily be combined with the method for weighted
nilpotent presentations. This promising approach yields a significant speed-up
as the following table shows. The timings were obtained on an Intel
Pentium 4 processor with a clock-speed of $2.4$~GHz. 
\begin{center}
  \begin{tabular}{cccccc}
\toprule
     Quotient & \verb|#gens| &\verb|Usual| & \verb|Solv| & \verb|Weight| & \verb|Solv+Weight| \\
\midrule
$\mbox{BSV}$, class $25$ & 106 & 0:00:05 & 0:00:04 & 0:00:01 & 0:00:01\\
$\mbox{BSV}$, class $35$ & 179 & 0:01:35 & 0:01:06 & 0:01:06 & 0:00:48\\
$\mbox{BSV}$, class $40$ & 219 & 0:04:26 & 0:03:00 & 0:03:22 & 0:02:25\\
$\mbox{BSV}$, class $45$ & 259 & 0:10:27 & 0:06:54 & 0:08:28 & 0:06:05\\
$\mbox{BSV}$, class $50$ & 301 & 6:31:17 & 3:52:36 & 6:30:13 & 4:43:52\\
\midrule
$\Delta$, class ${35}$ & 185 & 0:00:31 & 0:00:31 & 0:00:02 & 0:00:02\\
$\Delta$, class ${80}$ & 609 & 1:19:22 & 1:15:03 & 0:29:48 & 0:27:36\\
$\Delta$, class $100$  & 821 & 8:25:37 & 7:39:54 & 5:45:40 & 5:18:08\\
\bottomrule
  \end{tabular}
\end{center}

The method \verb|Usual| denotes the algorithm in [14, p. 424]
for polycyclic presentations, the method
\verb|Solv| denotes our new method, the method \verb|Weight|
denotes the method for weighted nilpotent presentation as 
in [14, p. 431], and the method \verb|Solv+Weight| denotes
the combination of both of the latter methods. The number \verb|#gens|
denotes the number of generators of the considered polycyclic presentation. 
In summary, our method always yields here a significant
speed-up compared with the standard method for polycyclic groups. \\ \\
{\it Acknowledgement}. We thank Michael Vaughan-Lee for many useful comments
and for having provided us with a simpler proof for Lemma 1.

\end{document}